\documentclass[11pt]{article}
\usepackage{graphicx}
 \usepackage[reqno, namelimits, sumlimits]{amsmath}
   \usepackage{amsfonts}
   \usepackage{amssymb}
  \usepackage{amsthm}

\newtheorem{thm}{Theorem}[section]

\newtheorem{remark}[thm]{Remark}

\def\R{{\mathbb R}}

\font\textmsbm=msbm10
\newfam\Bbbfam
\textfont\Bbbfam=\textmsbm
\def\bb{\fam\Bbbfam}\def\bb{\fam\Bbbfam}\def\bb{\fam\Bbbfam}
\date{ }
\begin{document}
\setcounter{page}{0}
\title{Some simple problems for the next generations.}

\maketitle

\medskip\author{ \center  Alain HARAUX \\ Sorbonne Universit\'es, UPMC Univ Paris 06, CNRS, UMR 7598, \\Laboratoire Jacques-Louis Lions, 4, place Jussieu 75005, Paris, France.\\
 haraux@ann.jussieu.fr\\}

\bigskip

\renewcommand{\abstractname} {\bf Abstract}
\begin{abstract}

A list of open problems on global behavior in time of some evolution systems, mainly governed by P.D.E,  is given together with some background information explaining the context in
which these problems appeared. The common characteristic of these problems is that they
appeared a long time ago in the personnal research of the author and received almost no answer
till then at the exception of very partial results which are listed to help the
readers' understanding of the difficulties involved.
\end{abstract}

\bigskip

\qquad {\small \bf AMS classification numbers:}  35B15, 35B40, 35L10, 37L05, 37L15
 
\vfill

\newpage

\noindent 
{\Large \bf  Introduction}

\bigskip

 \medskip \noindent
 
 Will the next generations go on studying mathematical problems? This in itself is an open question, but the growing importance of computer's applications in everyday's life together with the fundamental intrication of computer science, abstract mathematical logic and the developments of new mathematical methods makes the positive answer rather probable.

\bigskip

This text does not comply with the usual standards of mathematical papers for two reasons: it is a survey paper in which no new result will be presented and the results which we recall to motivate the open questions will be given without proof. 

\bigskip

It is not so easy to introduce an open question in a few lines. Giving the statement of the question is not enough, we must also justify why we consider the question important and explain why it could not be solved until now. Both points are delicate because the importance of a problem is always questionable and the difficulty somehow disappears when the problem is solved.

\bigskip

The questions presented here concern the theory of differential equations and mostly the case of PDE. They were encountered by the author during his research and some of them are already 40 years old. They might be considered purely academical by some of our colleagues more concerned by real world applications, but they are selected, among a much wider range of open questions, since their solution probably requires completely new approaches and will likely open the door towards a new mathematical landscape.

\section{Compactness and almost periodicity}

Throughout  this section, the terms ``maximal monotone operator" and ``almost periodic function" will be used without having been defined. Although both terms are by now rather well known, the definitions and main properties of these objects will be found respectively in the reference texts \cite{HB} and \cite{AP2}.
\\

One of my first fields of investigation was, in connection with the abstract  oscillation theory, the relationship between (pre-)compactness and asymptotic almost periodicity for the trajectories of an almost periodic contractive process. The case of autonomous processes (contraction semi-groups on a metric space) had been studied earlier in the Hilbert space framework by Dafermos and Slemrod \cite{D-S}, the underlying idea being that on the omega-limit set of a precompact trajectory, the semi-group becomes an isometry group. Then the situation resembles the simpler case of the isometry group generated on a Hilbert space $H$ by the equation 
$$ u'+Au(t) =  0 $$ with $$ A^* = -A $$
for which almost periodicity of precompact trajectories was known already from  L. Amerio quite a while ago (the case of vibrating membranes and vibrating plates with fixed bounded edge are special cases of this general result). 
\medskip

The case of a non-autonomous process, associated with a time-dependent evolution equation of the form $$ u'+A(t) u(t)\ni 0 $$ is not so good in general. In \cite {H1} I established an almost periodicity result for precompact trajectories of a {\it periodic} contraction process on a complete metric space, and in the same paper I exhibited a simple almost periodic (linear)  isometry process  on $\R^2$, generated by an equation of the form $$ u'+c(t)Ju(t) =  0 $$ with $J$ a $\pi\over2$- rotation around $0$,  for which no trajectory except $0$ is almost periodic. 

\medskip

Actually, while writing my thesis dissertation, I was specifically interested in the so-called ``quasi-autonomous" problem, and I met the following general question

\bigskip

{\bf Problem 1.1.} (1977)   Let A be a maximal monotone operator on a real
Hilbert space H, let $ f:{\R}\longrightarrow  H$ be almost periodic  and let $u$ be
a solution of  de $$ u'+A u(t)\ni f(t) $$ on $[0,
+\infty)$ with a  {\bf precompact range}. Can we conclude that $u$ is  asymptotically
almost periodic? 

\bigskip 

After studying a lot of particular cases in which the answer is positive ( $A = L$  linear, $A$  a subdifferential $\partial \Phi $ and some operators of the form $ L + \partial \Phi $ )  , I proved in \cite{HPP} that the answer is positive if $ H={\R}^N$  with 
$ N \le  2$. But the answer  is unknown  for  general maximal monotone operators even if $ H={\R}^3.$ 

\begin{remark} {\rm In  \cite{wrong} it is stated that the answer is positive for all $N$, but there is a  mistake in the proof, relying on a geometrical property which is not valid in higher dimensions, more specifically in 3D the intersection of the  (relative) interiors of two arbitrarily close isometric proper triangles can be  empty. Therefore the argument from \cite{HPP} cannot be used in the same way for $N\ge3$. }\end{remark}  

\begin{remark} {\rm The problem is also open even when $ A\in C^1(H, H)$, in which case the monotonicity just means $$\forall u\in H,\quad \forall v\in H, \quad (A'(u), v) \ge 0. $$ }\end{remark}  

\begin{remark} {\rm The answer is positive if $f$ is periodic, as a particular case of the main result of \cite {H1}. }\end{remark}  

\bigskip 

Since an almost periodic function has precompact range, studying the existence of almost periodic solutions requires some criteria for precompactness of bounded orbits. In the case of evolution PDE, precompactness is classically derived from higher regularity theory. For parabolic equations the smoothing effect provides some higher order regularity for $t>0$ for bounded semi-orbits defined on $\R^+$. In the hyperbolic case, although there is no smoothing effect in finite time, precompactness of orbits was derived by Amerio and Prouse \cite{AP} from higher regularity of the source and strong coercivity of the damping operator $g$ in the case of the semilinear hyperbolic problem
$$ u_{tt} -
\Delta  u +  g(u_t) = f(t, x)
\hbox { in }
{\R^+}\times
\Omega,\ \ \ \          u = 0    \hbox { on }    {\bb
R^+}\times
{\partial}\Omega $$ 
where $ \Omega $ be a bounded domain of ${\R}^N$.  But this method does not apply even in the simple case $g(v) = cv^3$ for $c>0, N\le 3$, a case where boundedness of all trajectories is known. The following question makes sense even when the source term is periodic in t and $g$ is globally Lipschitz continuous.

\bigskip

{\bf Problem 1.2.} (1978)   Let $ \Omega $ be a bounded domain
of
${\R}^N$ and $g$ a nonincreasing Lipschitz function. We consider the
semilinear hyperbolic problem
$$ u_{tt} -
\Delta  u +  g(u_t) = f(t, x)
\hbox { in }
{\R^+}\times
\Omega,\ \ \ \          u = 0    \hbox { on }    {\bb
R^+}\times
{\partial}\Omega $$ We assume that $ f:{\R}\longrightarrow  L^2(\Omega)$ is
continuous and periodic in t. Assuming  $$ u\in C_b({\R}^+,
H^1_0(\Omega))\cap C^1_b({\R}^+, L^2(\Omega)) $$
can we conclude that  $$ \bigcup_{t\geq 0} \{(u(t,.),u_t(t,.))\}\hbox{ is
precompact in }
  H^1_0(\Omega)\times L^2({\Omega})  ? $$
\begin{remark} {\rm The answer is positive in the following extreme cases

1) If g = 0 (by Browder-Petryshyn's theorem, there is a periodic solution, hence compact, and all the others are precompact by addition.)

2) If $g^{-1}$ is uniformly continuous, cf. \cite{Edin} , the result does not require Lipschitz continuity of $g$ and applies for instance to $g(v) = cv^3$ for $c>0, N\le 3$

\medskip

It would be tempting to ``interpolate", but even the case  $g(v) = v^+$ and $N=1$ already seems to be non-trivial. }\end{remark}

\begin{remark} {\rm The same question is of course also relevant when $f$ is almost periodic, and the result of  \cite{Edin} is true in this more general context. Moreover precompactness of bounded trajectories when $g = 0$ is also true when $f$ is almost periodic. This is related to a fundamental result of Amerio stating that if the primitive of an almost periodic function: $\R\longrightarrow H $  is bounded, it is also almost periodic. More precisely, if $ H $ is a Hilbert space and $L$ is a (possibly unbounded ) skew-adjoint linear operator with compact resolvent, let us consider a bounded solution (on $\R$ with values in $H$) of the equation 
$$ U'+ AU = F $$ where $F: \R\longrightarrow H$ is almost periodic . Then $ \exp (tA)U := V$ is a bounded solution of $$ V' = \exp (tA)F $$ and, since $\exp (tA)\psi $ is almost periodic as well as $\exp (- tA)\psi $ for any $\psi\in H$, by a density argument on generalized trigonometric  polynomials, it is immediate to check that a function $W: \R\longrightarrow H$ is almost periodic if and only if $\exp (tA)W: \R\longrightarrow H$  is almost periodic. Then Amerio's Theorem applied to $V$ gives the result, and this property applies in particular to the wave equation written as a system in the usual energy space. Then starting from a solution bounded on $\R^+$, a classical translation-(weak)compactness argument of Amerio gives a solution bounded on $\R$ of the same equation. We skip the details since this remark is mainly intended for experts in the field.} \end{remark}

\bigskip

\section{Oscillation theory}

Apart from the almost periodicity of solutions which provides a starting point to describe precisely the global time behavior of vibrating strings and membranes with fixed edge, it is natural to try a description of sign changes of the solutions  on some subset of the domain. Let us first consider  the basic equation  \begin{equation}\label{basic-eq}
u''+  Au(t)= 0,
\end{equation} where $V$ is a real Hilbert space, $A\in L(V, V')$ is a symmetric, positive, coercive operator and there is a second real Hilbert space $H$ for which $V\hookrightarrow H
=H' \hookrightarrow V'$ where the imbedding on the left is compact. In this case it is well known that all solutions $u\in C(\R , V) \cap C^1(\R , H)$ of \eqref{basic-eq}  are almost periodic : $\R\rightarrow V $  with mean-value $0$. Then for any form $\zeta\in V'$, the function $g(t): = \langle \zeta, u(t)\rangle$ is a real-valued continuous almost periodic function with mean-value $0$. It is then easy to show that either $g\equiv 0$, or there exists $M>0$ such that on each interval $J$ with $ |J | \ge M $, g takes both positive and negative values. We shall say that a number $M>0$ is a strong oscillation length for a numerical function $g\in L^1_{loc} (\R)$  if the following alternative holds: either $g(t) = 0$ almost everywhere, or for any interval $J$ with $ |J | \ge M $, we have 
$$ meas \{ t\in J, f(t) >0\} >0  \quad \hbox {and}\quad  meas \{ t\in J, f(t) <0\} >0.$$  As a consequence of the previous argument , under the above conditions on $H, V$ and $A$, for any solution $u\in C(\R , V) \cap C^1(\R , H)$ of \eqref{basic-eq} and for any $\zeta\in V'$, the function $g(t): = \langle \zeta, u(t)\rangle $ has some finite strong oscillation length $M = M(u, \zeta)$. 

\bigskip

In the  papers \cite {CH1, CH2, HK} the main objective was to obtain a strong oscillation length independent of the solution and the observation in various cases, including non-linear perturbations of equation \eqref{basic-eq}. A basic example is the vibrating string equation
\begin{equation}\label{NLwave}
 u_{tt} - u_{xx}+ g(t, u) = 0 \quad \hbox {in} \,\, \R\times (0, l), \quad u = 0 \,\, \hbox {on}\,\, \R\times \{0, l\}  
\end{equation} where $l>0$ and $g(t, .)$ is an odd non-decreasing function of $u$ for all $t$. Here the function spaces are $H= L^2(0, l)$ and $V= H^1_0(0, l)$.  Since any function of V is continous, a natural form $\zeta\in V'$ is the Dirac mass $\delta_{x_0}$ for some ${x_0}\in (0, l).$ It turns out that $2l$ is a strong oscillation length independent of the solution and the observation point ${x_0}$, exactly as in the special case $g= 0$, the ordinary vibrating string. Since in this case all solutions are $2l$-periodic with mean-value $0$  functions with values in $V$, it is clear that $2l$ is a strong oscillation length independent of the solution and the observation point ${x_0}$. The slightly more complicated $g(t, u) = au$ with $a>0$ is immediately more difficult since the general solution is no longer time-periodic, it is only almost periodic in $t$. The time-periodicity is too unstable and for an almost periodic function, the determination of strong oscillation lengths is not easy in general, as was exemplified in \cite {HK}. The oscillation result of \cite {CH1, CH2} is consequently not so immediate even in the linear case. In the nonlinear case,  it becomes even more interesting because the solutions are no longer known to be almost periodic.   \bigskip 

In dimensions $N\ge 2$, even the linear case becomes difficult. It has been established in \cite{HK} that even for analytic solutions of the usual wave equation in a rectangle, there is no uniform pointwise oscillation length common to all solutions at some points of the domain. One would imagine that it becomes true if the point is replaced by an open subset of the domain, but apparently nobody knows the answer to the exceedingly simpler following question: \medskip

{\bf Problem 2.1.} (1985)   Let $ \Omega = (0, 2l)\times (0, 2l)\subset
{\R}^2$. We consider the linear wave equation
$$ u_{tt} -
\Delta  u = 0
\hbox { in }
{\R}\times
\Omega,\ \ \ \          u = 0    \hbox { on }    {\bb
R}\times
{\partial}\Omega $$ Given $T>0$ , can we find a solution $u$ for which 
$$ \forall (t, x) \in [0, T]\times (0, l)\times (0, l),\quad  u(t, x) >0 ?  $$ Or does this become impossible for $T$ large enough? 

\bigskip 

Another simple looking intriguing question concerns the pointwise oscillation of solutions to semi linear beam equations, since the solutions of the corresponding linear problem oscillate at least as fast as those of the string equation:

\bigskip  {\bf Problem 2.2.} (1985)   \medskip We consider the
semilinear beam equation $$ u_{tt} + u_{xxxx} + g(u) = 0 
\hbox { dans } {\R}\times (0, 1) ,  \quad u =
u_{xx} = 0  \hbox { on  } {\R}\times \{0,1\}$$
with  g  odd and nonincreasing  with respect to  u. Is it possible for a solution
$u(t, .)$ to remain positive at some point $x_0$ on an arbitrarily long (possibly
unbounded) time interval?  

\bigskip Finally, let us mention a question on spatial oscillation of solutions to parabolic problems. Since the heat equation has a very strong smoothing effect on the data, and all solutions are analytic inside the domain for $t>0$, it seems natural to think that they do no accumulate oscillations and for instance in 1D, the zeroes of $u(t,.)$ will be isolated for $t>0. $  A very general result of this type, valid for semi linear problems as well  has been proved by Angenent \cite{Ang}. But as soon as $N\ge 2$, even the linear case is not quite understood. The answer to the following question seems to be unknown:

\bigskip

 {\bf Problem 2.3.} (1997)   Let $ \Omega \subset {\R}^N$
be a bounded open domain. We consider the heat equation
$$ u_{t} -
\Delta  u = 0
\hbox { in }
{\R}\times
\Omega,\ \ \ \          u = 0    \hbox { on }    {\bb
R}\times
{\partial}\Omega $$ For $t>0$ , we consider $$ {\cal E} = \{x\in \Omega,\quad  u(t, x)
\not= 0 \}$$  Is it true that $ {\cal E}$ has a finite number of connected
components? \medskip

\begin{remark} {\rm The solutions  $u$ of the {\it  elliptic }problem 
$$- \Delta  u + f(u)=0 \hbox { in } \Omega,\ \ \ \          u = 0    \hbox { on }   
{\partial}\Omega $$  are such that $ \{x\in \Omega,\quad  u( x)
\not= 0 \}$   has a finite number of connected components for a large class of functions $f$, cf. e.g. \cite{CHM}. Hence stationary solutions cannot provide a counterexample. } \end{remark}

\section{A semi-linear string equation}

There are in the Literature a lot of results on global behavior of solutions to Hamiltonian equations in finite and infinite dimensions. Apart from Poincar{\'e}'s recurrence theorem and the classical results of Liouville on quasi-periodicilty for most solutions of completely integrale finite dimensional hamiltonians, none of the recent results is easy and there is essentially nothing on PDE except in 1D. Even the case of semi linear string equations is not at all well-understood. While looking for almost-periodic solutions (trying to generalize the Rabinowitz theorem on non-trivial periodic solution) I realized that even precompactness of general solutions is unknown for the simplest semi linear string equation in the usual energy space :

\bigskip {\bf Problem 3.1.} (1976)   \medskip For the simple equation
  $$     u_{tt} - u_{xx} + u^3 = 0  \hbox { in } {\R} \times  (0,1) , \,\,
u = 0  \hbox { on }  {\R} \times  \{0,1 \}$$ the following simple
looking questions seem to be still open

\medskip  Question 1.    Are there solutions which converge weakly to 0 as
time goes to infinity?

\medskip  Question 2.    If $(u(0,.), u_t(0, .))\in H^2((0, 1)) \cap
H_0^1((0, 1))\times H_0^1((0, 1)): = {\cal V} $, does $(u(t,.), u_t(t, .))$ remain
bounded in
${\cal V} $ for all times?

\bigskip  \begin{remark} {\rm To understand the difficulty of the problem, let us just mention that the equation 
 $$     iu_{t} + |u|^2 u = 0  \hbox { in } {\R} \times  (0,1) ,\ \
u = 0  \hbox { on }  {\R} \times  \{0,1 \}$$ has many solutions tending weakly to $0$ and, although the calculations are less obvious, the same thing probably happens to $$     u_{tt} + u^3 = 0  \hbox { in } {\R} \times  (0,1) , \,\,
u = 0  \hbox { on }  {\R} \times  \{0,1 \}$$ Hence the problem appears as a competition between the ``good" behavior of the linear string equation and the bad behavior of the distributed ODE associated to the cubic term. } \end{remark}

\bigskip  \begin{remark} {\rm If the answer to question 2 is negative, it means that, following the terminology of Bourgain \cite{Bourg}, the cubic wave equation on an interval is a weakly turbulent  system. Besides, weak convergence to $0$ might  correspond to an accumulation of  steep spatial oscillations of weak amplitude, not contradictory with the energy conservation of solutions.} \end{remark}

\bigskip  \begin{remark} {\rm In \cite{CHW1}-\cite{CHW3}, the authors investigated the problem 
\begin{equation}\label{intterm}u_{tt} - u_{xx}+ u\int_0^lu^2(t, x) dx = 0 \quad \hbox {in} \,\, \R\times (0, l), \quad u = 0 \,\, \hbox {on}\,\, \R\times \{0, l\}   \end {equation} which can be viewed as  a simplified model to understand the above equation. In this case, there is no solution tending weakly to $0$, and the answer to question 2 is positive. Interestingly enough, in this case the distributed ODE takes the form  $     u_{tt} + c^2(t) u= 0  $ , so that the solution has the form $a(x) u_1(t) + b(x) u_2(t)$ and remains in a two-dimensional vector space! This precludes both weak convergence to $0$ and weak turbulence.
} \end{remark}

\section {Rate of decay for damped wave equations}

Let us consider the  semilinear hyperbolic problem
$$ u_{tt} -
\Delta  u +  g(u_t) = 0\hbox { in }
{\R^+}\times
\Omega,\ \ \ \          u = 0    \hbox { on }    {\bb
R^+}\times
{\partial}\Omega $$ 
where $ \Omega $ be a bounded domain of ${\R}^N$ and $g$ is a nondecreasing function with $g(0)  =0$. Under some natural growth conditions on $g$, the initial value problem is well-posed and can be put in the framework of evolution equations generated by a maximal monotone operator in the energy space  $$ H^1_0(\Omega) \times L^2(\Omega)$$  

An immediate observation is the formal identity  $$ \frac {d}{dt} [\int_\Omega (u_t^2 + | \nabla u|^2 ) dx] = - 2\int_\Omega g(u_t)u_t dx \le 0$$ showing that the energy of the solution is non-increasing. When $g(s) = cs$ with $c>0$, one can prove the exponential decay of the energy by a simple calculation involving a modified energy function 
$$ E_\varepsilon (t) = \int_\Omega (u_t^2 + |\nabla u|^2 ) dx + \varepsilon \int_\Omega uu_t dx $$  The exponential decay is of course optimal since $$ \frac {d}{dt} [\int_\Omega (u_t^2 + | \nabla u|^2 ) dx] = - 2\int_\Omega cu_t ^2dx \ge -2c \int_\Omega (u_t^2 + | \nabla u|^2 ) dx$$

A similar calculation  can be performed if $ 0< c\le g'(s) \le C $ , and the result is even still valid for 
$g(s) = cs + a |s|^\alpha s $ under a restriction on $\alpha>0 $ depending on the dimension. 

\bigskip 

More difficult, and somehow more interesting, is the case 

$$g(s) =  a |s|^\alpha s, \quad a>0, \alpha>0 $$
in which under a restriction relating $\alpha$ and $N$, various authors (cf. e.g. \cite{Nakao}, \cite{HZ} and the references therein) obtained the energy estimate 

$$ \int_\Omega (u_t^2 + |\nabla u|^2 ) dx \le C (1+t) ^{- \frac{2}{\alpha}}$$  But now the energy identity only gives $$ \frac {d}{dt} [\int_\Omega (u_t^2 + | \nabla u|^2 ) dx] = - 2\int_\Omega a|u_t |^{\alpha + 2}dx $$ while to prove the optimality of the decay we would need something like $$ \frac {d}{dt} [\int_\Omega (u_t^2 + | \nabla u|^2 ) dx] \ge - C( \int_\Omega u_t ^{ 2}dx) ^{1 +\frac{\alpha}{2} } $$ Unfortunately the norm of $ u_t$ in $L^{\alpha + 2}$ cannot be controlled in terms of the $L^2$ norm, even if strong restrictions on $u_t$ are known. If $u_t$ is known to be bounded in a strong norm, let us say an $L^p$ norm with $p$ large, we can derive a lower estimate of the type 
$$ [\int_\Omega (u_t^2 + | \nabla u|^2 ) dx]\ge \delta (1+t) ^{- \beta} $$ for some $\beta > \frac{2}{\alpha}$. But even $p= \infty$ does not allow to reach the right exponent. 

\bigskip In 1994, using special Liapunov functions only valid for $N = 1$, the author ( cf. \cite{Below}) showed that for all sufficiently regular non-trivial initial data, we have the estimate
$$ \int_\Omega (u_t^2 + |\nabla u|^2 ) dx \ge C (1+t) ^{- \frac{3}{\alpha}}$$
In general, for $N>2$ , some estimate of the form 

$$ \int_\Omega (u_t^2 + |\nabla u|^2 ) dx \ge C (1+t) ^{- K}$$
will be obtained if the initial data belong to $D(-\Delta)\times H^1_0(\Omega)$ and 
$\alpha < \frac{4}{N-2}$ . But we shall have in all cases $K > \frac{4}{\alpha}$ and $K$ tends to infinity when $\alpha $ approaches the value $\frac{4}{N-2}$ . 

\bigskip

\bigskip  \begin{remark} {\rm It is perfectly clear that none of the above partial results is satisfactory, since for analogous systems in finite dimensions, of the type 

$$u '' + Au + g(u')$$  with $A$ symmetric, coercive ,  $ (g(v), v) \ge c |v|^{\alpha + 2}$ and $ |g(v)|\le C |v|^{\alpha + 1}$, the exact asymptotics of any non-trivial solution is 
$$ |u'|^2 + |u|^2\sim  (1+t) ^{- \frac{2}{\alpha}}$$
Moreover, an optimality result of the decay estimate has been obtained in 1D  by P. Martinez and J. Vancostenoble \cite{Van-Mart} in the case of a boundary damping for which the same upper estimate holds. The difference is that inside the domain, an explicit formula  gives a lot of information on the solution. }\end{remark}

\bigskip {\bf Problem 4.1.}   \medskip For the  equation $$ u_{tt} -
\Delta  u +  g(u_t) = 0\hbox { in }
{\R^+}\times
\Omega,\ \ \ \          u = 0    \hbox { on }    {\bb
R^+}\times
{\partial}\Omega $$ with $$g(s) =  a |s|^\alpha s, \quad a>0, \alpha>0 $$
\medskip  Question 1.    Can we find a solution $u$  for which $$ |\int_\Omega (u_t^2 + |\nabla u|^2 ) dx \sim  (1+t) ^{- \frac{2}{\alpha}}? $$
\medskip  Question 2.    Can we find a solution $u$  for which the above property is {\it not} satisfied?

\bigskip  \begin{remark} {\rm Both questions seem to be still open for any domain and any $\alpha>0$. 
} \end{remark}

\section {The resonance problem for  damped wave equations with source term}

To close this short list, we consider the semilinear hyperbolic problem with source term
$$ u_{tt} -
\Delta  u +  g(u_t) = f(t, x)
\hbox { in }
{\R^+}\times
\Omega,\ \ \ \          u = 0    \hbox { on }    {\bb
R^+}\times
{\partial}\Omega $$ 
where $ \Omega $ be a bounded domain of ${\R}^N$. 
We assume that the exterior force $f(t, x)$ is bounded with values in $L^2(\Omega)$, In this case, all solutions $U = (u, u_t)$ are locally bounded on $(0, T)$ with values in the energy space $ H^1_0(\Omega) \times L^2(\Omega)$. The question is what happens as $t $ tends to infinity. 

\bigskip 

When $g(s) $ behaves like a super linear power $ |s|^\alpha s$  for large values of the velocity, it follows from a method introduced by G. Prouse \cite{Prouse} and extended successively by many authors, among which M. Biroli  \cite {Biroli}, \cite {Biroli-H} and the author of this survey, that the energy of any weak solution remains bounded for $t$ large, under the restriction $ \alpha (N-2) \le 4$ . Then many attempts  were tried to avoid this growth assumption. Many partial results were obtained under additional conditions ($ f $ bounded in stronger norms, $f $ anti-periodic, higher growths for $N\le 2$, cf e.g. \cite{Manus}, \cite{Manus2}, \cite{TG} ). But the following basic question remains open: 

\bigskip {\bf Problem 5.1.}   \medskip Assume $N\ge 3$,  $$g(s) =  a |s|^\alpha s, \quad a>0, \,\,\alpha>\frac {4}{N-2} $$
Is it still true that the energy of all solutions remains bounded for any exterior force $f(t, x)$  bounded with values in $L^2(\Omega)$? 

\bigskip  \begin{remark} {\rm The positive boundedness results require a weaker boundedness condition on $f$, it is sufficient that it belongs to a Stepanov space $S^p(\R, L^2(\Omega)$ with $p>1$. The first results in the direction were actually published by G. Prodi in 1956, so that the problem is about  60 years old...

} \end{remark}

\providecommand{\bysame}{\leavevmode\hbox to3em{\hrulefill}\thinspace}

\bibliographystyle{amsplain}


\providecommand{\bysame}{\leavevmode\hbox to3em{\hrulefill}\thinspace}
\providecommand{\MR}{\relax\ifhmode\unskip\space\fi MR }
\providecommand{\MRhref}[2]{%
  \href{http://www.ams.org/mathscinet-getitem?mr=#1}{#2}
}
\providecommand{\href}[2]{#2}
\begin{thebibliography}{}

\end{thebibliography}


\begin{thebibliography}{10}

\bibitem{AP} L. Amerio and G. Prouse, \emph{Uniqueness and almost-periodicity theorems for a non linear wave equation}, Atti Accad. Naz. Lincei Rend. Cl. Sci. Fis. Mat. Natur. \textbf{8}, 46 (1969) 1--8.

\bibitem{AP2} L. Amerio and G. Prouse, \emph{Almost-periodic functions and functional equations}, Van Nostrand Reinhold Co., New York-Toronto, Ont.-Melbourne 1971 viii+184 pp. 

\bibitem{Ang}  S. Angenent,  \emph{The zero set of a solution of a parabolic equation}, J. Reine Angew. Math. \textbf{390} (1988), 79--96. 

\bibitem{Biroli}  M. Biroli  ,  \emph{Bounded or almost periodic solution of the non linear vibrating membrane equation},  Ricerche Mat.  \textbf{22} (1973), 190--202. 

\bibitem{Biroli-H}  M. Biroli  and A. Haraux,  \emph{Asymptotic behavior for an almost periodic, strongly dissipative wave equation},  J. Differential Equations \textbf{38} (1980), no. 3, 422--440.

\bibitem{BVN} S. Bochner and J. Von Neumann, \emph{On compact solutions of operational-differential equations}, Ann. of Math. \textbf{2}, 36 (1935), 255--291.


\bibitem{Bourg} Jean Bourgain, \emph{On the growth in time of higher Sobolev norms of smooth solutions of Hamiltonian }, Internat. Math. Res. Notices \textbf{ 6}  (1996), 277--304. 

\bibitem{HB} H. Brezis,  \emph{Op{\' e}rateurs maximaux monotones et semi-groupes de contractions dans les espaces de Hilbert} (French) North-Holland Mathematics Studies, No. 5., Amsterdam-London; American Elsevier Publishing Co., Inc., New York, 1973. vi+183 pp. 

\bibitem{CH1} T. Cazenave and A. Haraux, \emph{Propri{\' e}t{\' e}s oscillatoires des solutions de certaines {\' e}quations des ondes semi-lin{\' e}aires}, C. R. Acad. Sci. Paris Ser. I Math. \textbf{298} (1984), no. 18, 449--452.

\bibitem{CH2} T. Cazenave and A. Haraux, \emph{Oscillatory phenomena associated to semilinear wave equations in one spatial dimension}, Trans. Amer. Math. Soc. \textbf{300} (1987), no. 1, 207--233.

\bibitem{CHW1} T. Cazenave, A. Haraux and F.B. Weissler,  \emph{Une {\' e}quation des ondes compl{\` e}tement int{\' e}grable avec non-lin{\' e}arit{\' e} homog{\` e}ne de degr{\' e} trois}, C. R. Acad. Sci. Paris Ser. I Math. \textbf{313} (1991), no. 5, 237--241.

\bibitem{CHW2} T. Cazenave, A. Haraux and F.B. Weissler, \emph{A class of nonlinear, completely integrable abstract wave equations}, J. Dynam. Differential Equations \textbf{5} (1993), no. 1, 129--154.

\bibitem{CHW3} T. Cazenave, A. Haraux and F.B. Weissler,  \emph{Detailed asymptotics for a convex Hamiltonian system with two degrees of freedom.}, J. Dynam. Differential Equations \textbf{5} (1993), no. 1, 155--187.

\bibitem{CHM}  M. Comte,  A. Haraux and P.  Mironescu, \emph{Multiplicity and stability topics in semilinear parabolic equations}, Differential Integral Equations \textbf {13} (2000), no. 7-9, 801--811.

\bibitem{D-S} C.M. Dafermos and M.  Slemrod, \emph{Asymptotic behavior of nonlinear contraction semigroups},  J. Functional Analysis \textbf {13 }(1973), 97--106.

\bibitem{TG} T. Gallouet, \emph{Sur les injections entre espaces de Sobolev et espaces d'Orlicz et application au comportement {\`a} l'infini pour des {\' e}quations des ondes semi-linéaires} 
(French), Portugal. Math. \textbf{42} (1983/84), no. 1, 97-- 112 (1985). 

\bibitem {H1} A. Haraux,  \emph{Asymptotic behavior of trajectories for some nonautonomous, almost periodic processes}, {\it J. Diff. Eq. } {\bf 49}(1983), no. 3, 473--483.

\bibitem {H2} A. Haraux,  \emph{A simple almost-periodicity criterion and applications}, {\it J. Diff. Eq. } {\bf 66} (1987), no. 1, 51--61. 

\bibitem {HPP} A. Haraux, \emph{Asymptotic behavior for two-dimensional, quasi-autonomous, almost periodic evolution equations}, {\it J. Diff. Eq. } {\bf 66} (1987), no. 1, 62--70.

\bibitem {Below} A. Haraux, \emph{ $ L^p$ estimates of solutions to some non-linear wave equations in one space dimension}, Int.J. Math. Modelling and Numerical Optimization {\bf 1} (2009), Nos 1-2,  p. 146-- 154.

\bibitem {OscPort} A. Haraux, \emph{On the strong oscillatory behavior of all solutions to some second order evolution equations}, Port. Math. {\bf 72} (2015), no. 2, 193--206. 

\bibitem {Edin} A. Haraux,  \emph{Almost-periodic forcing for a wave equation with a nonlinear, local damping term},  Proc. Roy. Soc. Edinburgh Sect. A {\bf 94} (1983), no. 3--4, 195--212.

\bibitem {Manus} A. Haraux, \emph{Nonresonance for a strongly dissipative wave equation in higher dimensions},  Manuscripta Math. {\bf 53 }(1985), no. 1-2, 145--166. 

\bibitem {Manus2} A. Haraux, \emph{ Anti-periodic solutions of some nonlinear evolution equations}, Manuscripta Math. 63 (1989), no. 4, 479--505.  

\bibitem {Dieud} A. Haraux,  \emph{Semi-linear hyperbolic problems in bounded domains},  Math. Rep. 3 (1987), no. 1, i--xxiv and 1--281. 

 \bibitem{HK} A. Haraux  and V. Komornik, \emph{ Oscillations of anharmonic Fourier series and the wave equation},  Rev. Mat. Iberoamericana \textbf{1} (1985), no. 4, 57--77.
 
 \bibitem{HZ} A. Haraux and E. Zuazua,\emph{ Decay estimates for some semilinear damped hyperbolic problems}, Arch. Rational Mech. Anal. {\bf 100} (1988), no. 2, 191--206. 
 
 \bibitem {wrong} Z. Hu  and A.B.  Mingarelli,  \emph{Almost periodicity of solutions for almost periodic evolution equations}, {\it Differential Integral Equations} {\bf 18} (2005), no. 4, 469--480. 

\bibitem{Mucken} C. F. Muckenhoupt, \emph{Almost periodic functions and vibrating systems}, Journal of Mathematical Physics \textbf{8 } (1928--1929), 163--199.


\bibitem{Nakao} M.Nakao, \emph{Asymptotic stability of the bounded or almost periodic solution of the wave equation with nonlinear dissipative term}, J. Math. Anal. Appl. {\bf 58} (1977), no. 2, 336--343.

\bibitem{Prod} G. Prodi, \emph{Soluzioni periodiche di equazioni a derivate parziali di tipo iperbolico non lineari} (Italian),  Ann. Mat. Pura Appl. (4) {\bf 42} (1956), 25--49.

\bibitem{Prouse}  G. Prouse, \emph{Soluzioni limitate dell'equazione delle onde non omogenea con termine dissipativo quadratico } (Italian),  Ricerche Mat. {\bf 14} (1965),  41--48. 

\bibitem{Van-Mart} J. Vancostenoble and P. Martinez,  \emph{ Optimality of energy estimates for the wave equation with nonlinear boundary velocity feedbacks}, SIAM J. Control Optim. 39 (2000), no. 3, 776--797 (electronic). 

\bibitem{Z}  S. Zaidman,  \emph{ Solutions presque-p{\'e}riodiques des {\'e}quations diff{\'e}rentielles abstraites}(French), Enseign. Math.  (2) {\bf 24} (1978), no. 1-2, 87--110.


\


\end{thebibliography}

\providecommand{\bysame}{\leavevmode\hbox to3em{\hrulefill}\thinspace}
\providecommand{\MR}{\relax\ifhmode\unskip\space\fi MR }
\providecommand{\MRhref}[2]{%
  \href{http://www.ams.org/mathscinet-getitem?mr=#1}{#2}
}
\providecommand{\href}[2]{#2}

\end{document}